\documentclass[12pt,letterpaper]{article}
\usepackage{t1enc}
\usepackage[latin1]{inputenc}
\usepackage[english]{babel}
\begin{document}
\date{}
\title{The  Period Function of Second Order Differential Equations }
\author{A.  Raouf  Chouikha \footnote
{Universite Paris 13 LAGA UMR 7539 Villetaneuse 93430,  e-mail: chouikha@math.univ-paris13.fr}
}
\maketitle

\begin{abstract}
 We  interest in the behaviour of the period function for  equations of the type $u'' + g(u) = 0$ and $u'' + f(u)u' + g(u) = 0$ with a center at the origin $0$. $g$ is a function of class $C^k$.  For the conservative case, if $k \geq 2$ one shows that the Opial criterion is the better one among those for which these the necessary condition $g''(0) = 0$ holds.  In the case where $f$ is of class $C^1$ and $k \geq 3$ , the Lienard equations  \ $ u'' + f(u) u' + g(u) = 0$ \ may have a monotonic period function if $g'(0) g^{(3)}(0) - \frac{5}{3} {g''}^{2}(0) - \frac{2}{3} {f'}^{2}(0) g'(0) \neq 0$ in a neighborhood of $0$.
\\ {\it Key Words and phrases:} \ period function, monotonicity, isochronicity, Lienard equation, polynomial systems.\footnote
{1991 Mathematics Subject Classification  \ 34C25, 34C35}
 
\end{abstract}
\maketitle
\vspace{1cm}

\section{Introduction}
The general question of monotonicity of the period function has been treated by several authors in the special case when the phase flow has a Hamiltonian form. In particular, there is the case when the differential equation is of conservative Newtonian type $$(E_c) \hspace{3cm} x'' + g(x) = 0.$$ For a review of the different criteria, one will be able to consult the  papers [R], [C], [Sc], [Ch-C] and the references herein.\\
C. Chicone [C2] points out a special interest of the behavior of the period function in understanding the Neumann boundary value problem \ $ x'(0) = 0, \ x'(T) = 0$.\ A solution of this problem corresponds to one of the periodic orbits in the phase plane that starts at time $t = 0$ on the x-axis and returns to the x-axis after time $T$. The monotonicity ensures the uniqueness of the solutions with a fixed number of nodes $N$, which corresponds to such a periodic trajectory with minimum period $\frac{2T}{N}$. This situation may occur for autonomous systems other than conservative ones. The interval of values of the period determines the interval of initial conditions for which the equation has a solution.\\

In outsides of this position, we know only little of cases, where the behavior of the period function perfectly is known as well as the isochronicity of the center.\\
The Lienard equation $$(L) \hspace{2cm} x'' + f(x) x' + g(x) = 0 ,$$ is not in Hamiltonian form and new result on the monotonicity of their period function can appear particularly interesting. \\

 M. Sabatini [S] proved if the functions \ $f$\ and \ $g$\ are analytic odd functions such that $f(0) = g(0) = 0, g'(0) > 0,$\ the origin $0$ is an isochronous center if and only if $$g(x) = g'(0) x + \frac{1  }{ x^3} [\int_0^x s f(s) ds]^2$$
These hypotheses in particular impose the strong condition \ $g''(0) = 0$.\ This one will be the central point of our discussion, where one will attempt to overcome this hypothesis.\\

Concerning the above problem, the interesting question that remains  is the one to know if the only isochronous centers  for  Lienard systems are those given by [S], as believe certain authors.  For the moment, one not know  examples of other isochronous systems.\\

This paper is organized as follows. In the second part we examine
monotonicity conditions of the period function for the Lienard
equation where the functions $f$ and $g$ are such that $f \in C^1(J,R)$ and $g \in C^3(J,R)$ on the interval $J = (a,b),\quad a < 0 < b$\ and 
satisfying some conditions so that the origin $0$ is a center of
this equation. It implies in particular \ $g'(0) > 0.$\   The period \ $T \equiv T(\gamma )$ \ where \ $\gamma $\ is a cycle surrounding \ $0$\ is defined.\\

We prove that  
$$g'(0) g^{(3)}(x) - \frac{5}{3} {g''}^2(x) - \frac{2}{3} {f'}^2(x) g'(0) \neq 0$$
implies the monotonicity of the period function \ $T$ \ in a neighborhood of $0$ . 
In particular,   $g'(0) g^{(3)}(0) - \frac{5}{3} {g''}^2(0) - \frac{2}{3} {f'}^2(0) = 0$ \ and \ $f'(0) g''(0) - g'(0) f''(0) = 0$  \ are necessary conditions for the  center $0$ to be isochronous. \\ When $g$ is of class $C^3$,  our criteria  appears in a certain direction to be more general than those given
by Sabatini [S] for the period function to be monotonic.\\ In particular, if  $f \equiv 0$, 
 we find the condition  $g'(0) g^{(3)}(0) - \frac{5}{3} {g''}^2(0)  \neq 0$ 
which implies the monotonicity of the period for the conservative case. This  agrees with Schaaf condition, [R] $$g'(x) g^{(3)}(x) - \frac{5}{3} {g''}^2(x) \neq 0$$ in  a neighborhood of $0$. \\ 

As a corollary of this result, we obtain the following\\

{\bf Theorem} \qquad {\it Let the origin $0$ be a center of  equations \ $x'' + g(x) = 0$\ and \ $x'' \pm  f(x) x' + g(x) = 0$\ 
Supppose the two functions  $f  \in C^1([a,b]),$  and $g \in C^3([a,b]),$  verify \ $f(0) = g(0)  = 0, g'(0) > 0.$    \\ If the period function \ $T$ \ of \ $x'' + g(x) = 0 $
\ is increasing in a neighborhood of $0$,  then the period function of \ $x'' \pm  f(x) x' + g(x) = 0$\ is also increasing in a neighborhood of $0$.} \\

As a consequence me may deduce that under the hypothesis \ $g''(0)\neq 0$ \ a necessary condition for the Lienard sytem to have an  isochronous center, is the period of the conservative associated system has to be decreasing.\\

We also notice that under restrictive hypotheses, the function $g(x)$ is of class $C^k, \ k \geq 3$ and $f \in C^2$, it is possible to extend some results of [S]. For that, one has to suppose $$g''(0) = 0 , \quad f''(0) = 0 \quad {\it and}\quad g'(x) g^{(3)}(x) - \frac{2}{3} {f'}^2(x) \neq 0$$ in a neighborhood of $0$, for the origin to be a center and the period function to be monotonic. 

Nevertheless, our result presents the advantage of to be more
natural. In this direction that if $g''(0) = 0$ it reduces himself to the one of
Z. Opial [O] when $f \equiv 0$.\\
 This  point will be studied early in the first part, where one will remark that in the conservative case the criterion given by F. Rothe appears to be  the better one since it does not demand the above condition $g''(0) = 0$.\\
  \\

\section {On the monotonicity conditions of the period}
In this section, we give other results concerning the monotonicity of the period energy-function 
of the conservative system. In particular, we prove that if $g$ is of class $C^ 2$, we cannot improve the Opial result without  supposing $g''(0) \neq 0$ . If  $g$ is of class $C^ 3$ and to obtain the same result  we have to suppose  additional conditions.\\ Moreover, we prove that Opial result is the better one 
among those for which these  $g''(0) = 0$. \\
Let us consider the conservative system 

 \begin{eqnarray}
 x'' + g(x) = 0.
 \end{eqnarray}

$0$ is a center for (2.1), $g \in C^1(J,R)$, where $J$ is a real
interval containing $0$ and satisfying  $g(0) = 0$ and $g'(0) >
0$. Let $G$ the primitive of $g$, verifying :
$G(0)=0 \qquad {\it and} \qquad  G(a)=G(b)=c, $ with $\ \alpha ,\beta $ satisfy the
inequalities\quad $\alpha <a<0\quad $ and\quad $0<b<\beta .$ \
Thus, there exists a constant\quad $\delta ,\quad $such that \quad
$ 0<c<\delta .$ \\ The energy-period function is
$T(r)=\sqrt{2}\int_a^b\frac{dx}{\sqrt{G(r)-G(x)}}$, \quad where
$G(r) = c$.\\

We prove the following , which slightly improves Opial result ([O], Theorem 8) \\

 {\bf  Proposition 1.}\qquad
 {\it Let a function \ $g$,\ and $J$  an interval containing $0$ such that  $g(0) = 0$ . Suppose  the function $g$ satisfies   \quad $x g(x) > 0$ for $x \in J$, \ then for each following condition \\
(i) - \ $g \in C^1(J,R)$,\ and \ $g'(0) > 0$.\\
(ii) - \ $g \in C^2(J,R)$,\ and \ $g'(0) > 0$\ and \ $g''(0) = 0$,\\ we have

\begin{eqnarray}
{\it if} \  \frac{d }{dx} (\frac{g }{x}) \ {\it is}
\cases{
{\it \ positive \ for} \   x \in J \  {\it and} \  x < 0, & \cr
{\it \ negative \ for} \  x \in J \ {\it and} \ x > 0, & \cr}
\end{eqnarray}

then the period function  $T$  of (2.3) is strictly increasing in a neighborhood of $0$;

\begin{eqnarray}
{\it if} \  \frac{d }{dx} (\frac{g }{x}) \ {\it is} 
\cases{
{\it \ positive \ for} \ x\in J \ {\it and} \ x > 0, & \cr
{\it \ negative \ for} x\in J \ {\it and} \ x < 0, & \cr}
\end{eqnarray}

then the period function $T$ of (2.3) is strictly decreasing in a neighborhood of $0$.  \\
Moreover, if \ $g \in C^3(J,R)$,\ then \ $ x g^{(3)}(x)$\ and \ $ \frac{d  }{dx}(\frac{g }{x})$ \ have a same sign in a neighborhood of $0$. }\\

{\bf Proof} \qquad Consider for $r\in (0,b)$, $$T_1(r) = \sqrt{2}\int_0^r\frac{dx}{\sqrt{G(r)-G(x)}}.$$The proves are similar if we start from hypothesis (2-2) or (2-3).
Suppose for example (2-2) is verified, let us prove that $T_1$ is strictly increasing on $(0,b)$. Following [O], consider $p = \frac{r'  }{r}$, where $0 < r < r' < b$.
 Define the new function $$h_p (x) = \frac{g(p x)  }{p}.$$ By (2-2),  for $x\in (0, \frac{b  }{p})$ we have $$g(x) = \frac{g(x)  }{x} x > \frac{g(px)  }{px} x = h_p (x).$$
  Let $I(x) = \int _0^x h_p (u) du $, we then have $G(r) - G(x) > I(r) - I(x).$
  We obtain $$T_1 (r) < \sqrt{2}\int_0^r\frac{dx}{\sqrt{I(r')-I(x)}} = \sqrt{2}\int_0^r\frac{dx}{\sqrt{G(r')-G(px)}} = T_1(r').$$ Then  $T_1$ is strictly increasing on $(0,b).$
  By  similarity , we prove that under hypothesis (2-2) $$T_2 (s) = \sqrt{2}\int_0^r\frac{dx}{\sqrt{G(s)-G(x)}}$$ is strictly decreasing on $(a,0)$.
   The period function $T$ is obviously $$T(r) = T_1(r) + T_2(s)$$ where $s\in (a,0)$ is such that $G(s) = G(r)$. \\ $T(r)$ is increasing because $G(x)$ is strictly increasing for $ x > 0$
   and strictly decreasing for $x < 0$.\\
Furthermore,
let the function \ $\psi(x) =  x (\frac{d}{dx}\frac{g}{x}))$,\ its
derivative is\ $\psi '(x) = g''(x) - {\frac{\psi (x)}{x}}$.\ By
the l'Hopital's rule we find $$2\ lim_{x\rightarrow 0}
(\frac{d}{dx}\frac{\psi }{x}) =  lim_{x\rightarrow 0} \psi '(x)
= g''(0) = 0.$$ Moreover, the second derivative of the function
$\psi (x)$ is $$\psi ''(x) = g^{(3)}(x) -  {\frac{\psi '(x)}{x}} +
{\frac{\psi (x)}{x^2}} = g^{(3)}(x) - {\frac{g'' (x)}{x}} + 2
{\frac{\psi (x)}{x^2}}.$$
 Thus, by l'Hopital's rule again we obtain  $$lim_{x\rightarrow 0} \psi ''(x) = 2 lim_{x\rightarrow 0} {\frac{\psi (x)}{x^2}} = lim_{x\rightarrow 0} {\frac{\psi'(x)}{x}} =  g^{(3)}(0)- lim_{x\rightarrow 0} {\frac{\psi (x)}{x^2}}.$$
So, $g^{(3)}(0)$ has the same sign as the function $\psi (x) =  x
(\frac{d}{dx}\frac{g}{x}))$ . \\ thus, $ x \frac{d  }{dx}(\frac{g }{x})$ has a constant sign only if we suppose \ $g''(0) = 0.$    

\bigskip

\subsection{Remark 1} 
{\bf  (i)} \quad In the convex case of the function\ $g(x)$  \ in the interval \ $(a,b)$,  \ Chow  and  Wang  [C-W] gave a  criterion of monotonicity for the period function . They prove that
 \begin{equation}
\label{C1} \left\{
\begin{array}{lll}
(i) & g^{\prime \prime }(x)>0, & {\it for} \qquad 
x\in (a,b), \\
(ii) & \Delta (x)=x(g^{\prime \prime }(0)g^{\prime }(x)-g^{\prime
}(0)g^{\prime \prime }(x))\geq 0 (\leq 0), & {\it for} \qquad x\in (a,b).
\end{array}
\right. (\mathcal{C}_1)
\end{equation}
implies

\begin{equation}
\label{C0} \left\{
\begin{array}{ll}
{H_0(x)={g(x)}^2+\frac{g^{\prime \prime }(0)}{3{g^{\prime }(0)}^2}{g(x)}
^3-2G(x)g^{\prime }(x) > 0 (< 0),} &  {\it for} \qquad x \in (a,b),  x \neq 0,
\end{array}
\right. (\mathcal{C}_0)
\end{equation}

which implies the monotonicity of the period function \ $T(r)$ for $0 < r < r_1$.\\

{\bf (ii)} \quad In the case where\ $g(x)$ \ is of class $C^2$ and non convex in the interval \ $(a,b)$,\ S.N. Chow and D. Wang prove a analogous result , where the function \ $\frac{d  }{dx}(\frac{g}{x})$\ in Proposition 1 is replaced by \ $x \frac{d^2  }{dx^2}g$ \  (their result is obviously weaker than Proposition 1 because $x \frac{d^2  }{dx^2}g$ and $\frac{d  }{dx}(\frac{g}{x})$ share the same sign). \\ But F. Rothe exhibited a stronger criterion, where this function is replaced by \\ $3{g''(0)}g'(x)^2 - {g''(0)}g(x)g''(x) - (3 {g'(0)^2  }) g''(x)$. The condition

\begin{equation}
\label{C4} \left\{
\begin{array}{ll}
H_{4}(x)=x[3{g''(0)}g^{\prime}(x)^{2}-{g''(0)}g(x)g^{\prime\prime}(x)-(3{g^{\prime
}(0)^{2}})g^{\prime\prime}(x)]\geq0 \ ,\\
{\it for} \qquad x\in(a,b),
\end{array}
\right. (\mathcal{C}_4)
\end{equation}

  denoted  $f_{4}$ in [R],  was proposed by
F. Rothe . \\ 
Remark that the
following assumption
$$
xg^{\prime\prime}(x)<0\quad  {\it for all }\quad x\neq 0
$$
implies also $(\mathcal{C}_{4})$ to hold. But, then necessarily implies $g^{\prime\prime}(0)=0,$ which is
a strong condition actually (see Proposition 2 below).\\

{\bf (iii)}\quad Recall that conditions $(\mathcal{C}_{0})$ and $(\mathcal{C}_{1})$ are
given by Chow and Wang. \\
The following criterion $(\mathcal{C}_{3})$ is due to R. Schaaf,
but \textit{(i)} is replaced by a weaker condition:

$$
{\it if} \qquad {g^{\prime}(x)=0} \ , \qquad {\it then} \qquad
{g(x)g^{\prime\prime}(x)<0.}$$

\begin{equation}
\label{C3} \left\{
\begin{array}{lll}
  H_3 (x)=5 {g^{\prime \prime }}^2(x) - 3g^{\prime }(x)g^{(3)
}(x) > 0 (< 0), \\  {\it for} \qquad {x\in (a,b), \ x \neq 0.}
\end{array}
\right. (\mathcal{C}_2)
\end{equation}

Note that it is equivalent to $(G^{\prime\prime}
)^{-2/3}$ being convex ($G$ is the primitive of $g$). \\

Notice that our following criterion $(\mathcal{C}_{5})$ (see [Ch-C]) which  is more inclusive than conditions $(\mathcal{C}_{1})$ and $(\mathcal{C}_{3})$ and with an additional assumption, more general than the above Rothe and Schaaf criteria (  respectively  $(\mathcal{C}_{4})$ and $(\mathcal{C}_{3})$ ) .

$$
\left\{
\begin{array}
[c]{lll}%
(i) & {g^{\prime\prime}(0)}[{3 {g^{\prime}(x)}}^{2}-g(x)g^{\prime\prime}(x)-3{g^{\prime
}(0)}^{2} g^{\prime\prime}(x) \neq 0 , &
{\it for} \qquad {{x\in({g^{\prime}}^{-1}(0),0),}} \\
(ii) & \frac{g^{\prime}(x) g^{\prime\prime}(0)} {g^{\prime\prime}
(x){g^{\prime}(0)}^{2}} \neq \frac{2G(x)}{{g(x)}^{2}}, & {\it for}\quad x\in(0,b), x \neq 0.
\end{array}
\right. {(\mathcal{C}_{5})}
$$

As we have seen in our preceding works [Ch-C], each of these conditions implies
that $(\mathcal{C}_{0})$ holds, which itself implies the monotonicity of the period.\\

Concerning the non convex case of the function $g(x)$ ( see remark (ii) above ), the following result gives relations between monotonicity conditions of the period function. But especially in the case where $g(x)$ is of class $C^2$ or $C^3$, it does to appear necessary conditions that are very restrictive.  What puts in evidence the role of the Rothe criterion.  Not only because it is the better one, but also because it does not require any supplementary conditions for the function \ $g(x)$. This result seems to not be known under this form, in any case it improves Theorem 1 of [R].\\

\subsection {\quad \bf Proposition  2} \qquad {\it Let \ $J$\ be an interval containig $0$. Suppose the function $g(x)$ is of class $C^2$  satisfying $x g(x) > 0$  for $x \in J$ and $g'(0) > 0$. Let the function $G(x)$ be the primitive of $g(x)$ verifying $G(0) = 0$ and $G(a) = G(b) ,\ a < 0 < b, \ a,b \in J$. Then we have the following implications\\
$$x g''(x) < 0 \quad {\it implies} \quad g^2(x) - 2G(x) g'(x) > 0 \quad {\it implies} \quad x (\frac{d}{dx}\frac{g}{x})) < 0.$$
Moreover, each of these conditions implies that the period function $T(r)$ of (2.3) is strictly increasing for $0 < r < r_1$.\\
$$x g''(x) > 0 \quad {\it implies} \quad g^2(x) - 2G(x) g'(x) < 0 \quad {\it implies} \quad x (\frac{d}{dx}\frac{g}{x})) > 0.$$
Moreover, each of these conditions implies that the period function $T(r)$ of (2.3) is strictly decreasing for $0 < r < r_1$.\\
A necessary condition to have any of these conditions is $g''(0) = 0$ .\\
 If we suppose in addition $g(x)$ is of class $C^3$, then $$g^{(3)}(x) < 0 (> 0) \ {\it and} \ g''(0) = 0 \ {\it implies} \ x g''(x) < 0 (> 0). $$}

\bigskip 

{\bf Proof} \qquad Notice that, since $G(x) > 0$ for $x \neq 0$ and  $$\frac{d}{dx}[g^2(x) - 2G(x) g'(x)] =  2 G(x) g''(x) ,$$
then \ $2x G(x) g''(x) < 0
\ ({\mbox resp.}> 0)$ \   implies \ $g^2(x) - 2G(x) g'(x) > 0
\ ({\mbox resp.}< 0)$.\  Because \ $x {\frac{d}{dx}[g^2(x) - 2G(x) g'(x)]} $ \ and \ $x g''(x)$\ have the same sign. So, we have proved the two first
implications.\\ 
The second implications \ $g^2(x) - 2G(x) g'(x) >
0 \ ({\mbox resp.}< 0)\ {\it implies} \ x
(\frac{d}{dx}\frac{g}{x}) < 0\  ({\mbox resp.}> 0)$ \ have been
proved by Rothe (see [R] Proposition 4 p. 138, according his
notations:  $h^+_3 \subset h^+_2$).\\

Remark that conditions  $$x g''(x) < 0 \ {\it resp.}> 0 \quad {\it and} \quad g^2(x) - 2G(x) g'(x) < 0 \ {\it resp.} > 0$$ independently imply \ $H_0(x)={g(x)}^2 + \frac{g^{\prime \prime }(0)}{3 {g^{\prime }(0)}^2} {g(x)}^3 -2G(x)g^{\prime }(x) \geq 0, \ ({\it resp.} \leq 0) $  and imply the
period function increasing (resp. decreasing), see [C-W] Corollary
2.3.\\
 On the other hand, condition \ $g^{(3)}(x) < 0 (> 0) \
{\mbox and} \ g''(0) = 0 $\ implies 

\begin{equation}
\begin{array}{lll}
  H_3 (x)=5 {g^{\prime \prime }}^2(x) - 3g^{\prime }(x)g^{(3)
}(x) > 0 (< 0), & {\it for} \qquad {x\in (a,b), x \neq 0.}
\end{array}
 (\mathcal{C}_3)
\end{equation}

holds, which itself implies $(\mathcal{C}_{0})$.\\ We proceed as in Proposition 1, 
let the function \ $\psi(x) =  x (\frac{d}{dx}\frac{g}{x}$,\  By
the l'Hopital's rule we find 
$$2\ lim_{x\rightarrow 0} (\frac{d}{dx}\frac{\psi }{x}) = g''(0) = 0.$$ We have seen that
$\psi (x)$ is $$\psi ''(x) = g^{(3)}(x) -  {\frac{\psi '(x)}{x}} +
{\frac{\psi (x)}{x^2}} = g^{(3)}(x) - {\frac{g'' (x)}{x}} + 2
{\frac{\psi (x)}{x^2}}.$$

 We have also calculated   $$lim_{x\rightarrow 0} \psi ''(x) =  g^{(3)}(0)- lim_{x\rightarrow 0} {\frac{\psi (x)}{x^2}}.$$
So, notice that $g^{(3)}(0)$ has the same sign as the function $\psi (x) =  x
(\frac{d}{dx}\frac{g}{x})$

Thus, $g^{(3)}(x) < 0 (> 0) \ {\it and} \ g''(0) = 0 \ {\it implies} \ x g''(x) < 0 (> 0), $ which implies \ $ x (\frac{d}{dx}\frac{g}{x})) < 0 (> 0)$.
 \\

From the latter remark, if \ $g^{(3)}(0) \neq 0$\ the function \ $x(\frac{d}{dx}\frac{g}{x})$\ has the same sign as \ $g^{(3)}(0),$\ we may deduce the following.

{\subsection {\bf Corrolary 1} \qquad {\it Let a function \ $g \in C^k (J,R), \ k \geq 2$\ where\  $J$\ is an interval containing $0$ such that \ $g(0) = g''(0) = 0, \ g^{(3)} (0) \neq 0$\ and  \ $g'(0) > 0$.\ Suppose in addition \ $x g(x) > 0 $\ for \ $x \in J$.\ Then we have  

\begin{eqnarray}
{\it if} \ g^{(3)}(0) < 0 \ {\it and }\ \frac{d  }{dx}(\frac{g }{x}) \ {\it is} \cases{ {\it \ positive \ for}\  x\in J \ {\rm  and}\  x < 0, &\cr
{\it \ negative \ for}\ x\in J \ {\it and}\ x > 0, &\cr  }
\end{eqnarray}
then the period function  $T$  of (2.1) is strictly increasing in a neighborhood of $0$;
 \begin{eqnarray}
{\it if} \ g^{(3)}(0) > 0 \ {\it and }\ \frac{d  }{dx}(\frac{g }{x}) \ {\it is} \cases{ {\it \ positive \ for} \ x\in J \ {\rm and} \ x > 0, &\cr
{\it \ negative \ for} x\in J \ {\rm and} \ x < 0, &\cr  }
\end{eqnarray}
then the period function $T$ of (2.1) is strictly decreasing in a neighborhood of $0$.  }\\

\vspace{0.5cm}

\section{Equation of Lienard type}
Let the equation
\begin{equation}
x'' + f(x) x' + g(x) = 0,
\end{equation}

such that $g(0) = 0$. Thus, $x \equiv 0$ is a trivial solution, and the origin is a singular point of the equivalent system
\begin{eqnarray}
   \cases{
   x' = - y  & \cr
   y' =  g(x) - f(x) y & \cr }
\end{eqnarray}

 Suppose $0$ is a center of (3.1), let $\gamma _0$ the {\it central region}, be the open connected set covered with cycles surrounded the center $0$. 
This periodic trajectories may be parametrized for exemple by choosing their initial values in the segment $(0, \pi)$ on the $x$ - axis.\\  Let  $T  : \gamma _0 \rightarrow R $ ,\  be the function defined  by associating to every point  $(x,0) \in \gamma _0$ the minimum period of the trajectory starting at $(x,0)$, to reach the negative x-axis.
 \ $T$ is the period function and is constant on cycles. We say $T$ is (strictly) increasing if, for every couple of cycles $\gamma _1$ and $\gamma _2 ,\  \gamma_1$ included in $ \gamma _2$ , we have \ $T(\gamma _1) \leq T(\gamma _2)$ \  ( $T(\gamma _1) < T(\gamma _2) $).\\ We say $0$ is isochronous center if $T$ is constant in a neighborhood of $0$.\\

Concerning this equation, different problems have been considered by several authors . In particular, existence, boundeness, uniqueness,  multiplicity of periodic solutions and related questions  were debated.\\
M. Sabatini [S]  interested in the monotonicity of the period function $T$ of a center $0$ of (3.1), or in the isochronicity of $0$.\\
Suppose the functions $f$ and $g$ are assumed to be $C^2$ functions on an interval $J$ which satisfy \ $f(0) = g(0) = 0,\ g'(0) > 0.$\ These assumptions ensure  the origin $0$ to be a center, so the period function $T$ is defined.\\ Since multiplication of the system of equations (3.2) by $\alpha ^{-1/2}$ does not change the nature (of monotonicity) of the period but only changes each period by a constant multiple. \\

More precisely, for any positive real number\ $\alpha $ \ equation (3.1) is equivalent to $$X'' + \frac{1  }{\sqrt{\alpha }} f(X) X' + \frac{1  }{\alpha }g(X) = 0$$ 
by the scaling 
$$x(t) = X ({\sqrt{\alpha }} t).$$ 

We are led to a system of the form  

\begin{eqnarray}
   \cases{
   x' = - \frac{1}{\alpha ^{-1/2}} y  & \cr
   y' =   \frac{1}{\alpha ^{-1/2}} g(x) -  \frac{1}{\alpha ^{-1/2}} f(x) y & \cr }
\end{eqnarray}

\bigskip

In fact it is equivalent to
another system more convenient to study,   We may show the following which agrees with Lemma 2 of
[S] for $g'(0) = 1$ \\

{\bf Lemma 1} \qquad {\it Suppose $f,g$ are continuous functions of
class $C^k, \ k \geq 1$ on $J$, an interval containing $0$ and $f(0) = 0$.
Let the function $C(x) = \frac{1  }{g'(0) }g(x)  - \frac{1  }{g'(0) x^3} [\int_0^x s f(s)
ds]^2$. Then the system
\begin{eqnarray}
\cases{
x' = y - \frac{1  }{x \sqrt{g'(0)}} \int_0^x s f(s) ds & \cr
y' = - g(x) - \frac{1  }{g'(0) {x}^3} [\int_0^x s f(s)
ds]^2 - \frac{y  }{x^2 \sqrt{g'(0)}} \int_0^x s f(s) ds & \cr}
\end{eqnarray}

  is of class $C^k, \ k \geq 1$ in a neighborhood of $0$ and equivalent to (3.1).}\\

{\bf Proof}\qquad Let us define\ $\psi (x) =  \frac{1  }{\sqrt{g'(0)} } \int_0^x s f(s) ds$.
\  By l'Hopital rule we get \\ $lim_{x\rightarrow 0} \frac{\psi
(x)  }{x^2} = \frac{f(0)  }{ {2 \sqrt{g'(0)} }}$. One also proves that the function
$ \frac{\psi (x)  }{x^2}$ is differentiable at $0$ and its
derivative takes the value $\frac{f'(0)  }{3{\sqrt{g'(0)} }}.$\\ Moreover, the
function $C(x)$ is obviously differentiable and we get $C'(0) = 1$. \\
A calculation gives $$C''(0) = \frac{g''(0) }{g'(0)} - \frac{2  }{3 g'(0)}f(0) f'(0) \quad {\it and} \quad C'''(0) = \frac{g'''(0) }{g'(0)} - \frac{2  }{3 g'(0)} f'(0)^2.$$

This proves the regularity of the system (3-3).\\
Furthermore,  let us consider $(x(t), y(t))$ a solution of (3-3),
then $x(t)$ is a solution of (3-1). Remark that for $x \neq 0$, we
get $$x [\frac{\psi (x)  }{x^2}]' = x^3 f(x) - \frac{2 x \psi (x)
}{x^3} = f(x) - 2\frac{\psi (x)  }{x^2}.$$ Moreover, by
differentiating $y = x' + x \frac{\psi (x)  }{x^2}$, we get $$x''
= -  \frac{1  }{g'(0)} g(x) -  \frac{1  }{\sqrt{g'(0)} }f(x) x'.$$
This is equivalent to (3.1) by scaling the time, here \ $\alpha = g'(0)$.\ This gives $$x( \frac{\tau }{\sqrt{g'(0)} }) = X(\tau)$$ we then obtain by differentiating with respect to $\tau$, and since $f$ and $g$ are independent on $t$, we have $f(x) = f(X)$ and $g(x) = g(X)$ 

$$X''= -   g(X) -  f(X) X'. $$

\bigskip

The following result specifies the behavior of the period function for the Lienard system  in the neighborhood of the center 0.  We have need nevertheless  hypotheses, to know \ $ f $ is $C^1$ and $g $ is $C^3$.  It allows of some to deduct several interesting corollaries.  Notably, simple conditions of monotonicity and to test in a quick way the isochronicity of the center.

\bigskip

\bf Theorem 1.} \qquad {\it Let a function $f  \in C^1([a,b])$ , and $g \in C^3([a,b])$ , verifying\ $f(0) = g(0)  = 0, g'(0) > 0$\ and  let the origin $0$ be a center of  $$x'' + f(x) x' + g(x) = 0.$$
  Then,
\begin{equation}
{\it if}\quad  g'(0) g^{(3)}(0) - \frac{5}{3} g''^2(0) - \frac{2}{3} f'^2(0) g'(0) \neq 0
\end{equation}

the period function $T$ of a periodic solution of this equation  is monotonic in a neighborhood of $0$.\\ More precisely,
$${\it if}\quad  g'(0) g^{(3)}(0) - \frac{5}{3} g''^2(0) - \frac{2}{3} f'^2(0) g'(0) <  0 $$ 
then $T$ is increasing in a neighborhood of $0$. 
$${\it if}\quad  g'(0) g^{(3)}(0) - \frac{5}{3} g''^2(0) - \frac{2}{3} f'^2(0) g'(0)  > 0   $$
then $T$ is decreasing in a neighborhood of $0$.} \\

\bigskip

{\bf Proof} \qquad To establish necessary conditions, we  give for that an expansion of the period function near the center. We will use  implicit function techniques.\\ Since the origin is a center,  orbits solutions starting on the positive x-axis must be  closed. Let \ $x(t,c), y(t,c)$ \ be a solution other than the origin of 

\begin{eqnarray}
 \cases{
x' = -   \sqrt{g'(0)} y  & \cr
y' =   \frac{1}{\sqrt{g'(0)}} g(x) -  \frac{1}{\sqrt{g'(0)}} f(x) y & \cr}
\end{eqnarray}

which take the value \ $x(0,c) = c , \ y(0,c) = 0 $\ at \ $t = 0.$\ Let us suppose $c$ is a positive constant closed to $0$. After a certain  time closed to \ $ \frac{2\pi}{\sqrt{g'(0)}}$ \ this solution will go around the origin and will again intersect the positive x-axis at \ $x(T,c)$.\ Consider the following functions dependant on $c$,  $$\phi (T,c) =  x(T,c) - c, \quad \psi (T,c) = y(T,c).$$
We will solve \ $\psi (T,c) = 0$\ for $T = T(c)$ \  a function on $c$ small. Thus,  \ $\phi $ is a function on $c$. Let  \ $\Phi (c) = \phi (T(c),c)$.\ We find that the position of return is \ $ x = c + \Phi (c).$ \ Thus,  the orbit is closed if and only if $$\Phi (c) = 0.$$
We will know the behavior of \ $\Phi (c)$ \  when $c$ tends to $0$, in calculating its first derivatives at $0$.\ At first, we have \ $\Phi (0) = 0,\  T(0) = \frac{2\pi}{\sqrt{g'(0)}}$ \ and the partial derivatives of $\phi $ and $\psi $ are $${\phi _T}(0,0) = x' (\frac{2\pi}{\sqrt{g'(0)}}, 0) = 0, \quad {\psi _T}(0,0) = y'(\frac{2\pi}{\sqrt{g'(0)}}, 0) = 0$$
$$  {\phi _c}(0,0) = {x_c}(\frac{2\pi}{\sqrt{g'(0)}}, 0) - 1, \quad {\psi _c}(0,0) = {y_c}(\frac{2\pi}{\sqrt{g'(0)}}, 0) .$$

The index $c$ or $T$ denote differentiation with respect to $c$ or $T$.\\
The derivatives $x_c(t,c)$ and $y_c(t,c)$ with respect to $t$ satisfy the Lienard system
\begin{eqnarray}
\cases {
x'_c = - {\sqrt{g'(0)}}y_c  & \cr
 y'_c =  \frac{1}{\sqrt{g'(0)}}g'(x)x_c - \frac{1}{\sqrt{g'(0)}}f(x) y_c - \frac{1}{\sqrt{g'(0)}}f'(x) x_c y & \cr }
\end{eqnarray}

with initial conditions \ $x_c(0,c) = 1, y_c(0,c) = 0$. According to our hypotheses, $g(0) = f(0) = 0$,  If we set  $c = 0$ the system becomes 

\begin{eqnarray}
\cases {
   x'_c(t,0) = - {\sqrt{g'(0)}}y_c  & \cr
   y'_c(t,0) =  {\sqrt{g'(0)}}x_c  & \cr }
\end{eqnarray}

It implies in particular, \ $x_c(t,0) = \cos({\sqrt{g'(0)}} t)$ \ and \ $y_c(t,0) = \sin({\sqrt{g'(0)}} t).$ \ Thus,\ $ \phi (0,0) = 0$ \ and \ $\psi (0,0) = 0$.\\
We now calculate the second derivatives $$\phi _{TT}(0,0) = x''(\frac{2\pi}{\sqrt{g'(0)}}, 0) = 0, \quad   \psi _{TT}(0,0) = y''(\frac{2\pi}{\sqrt{g'(0)}}, 0) = 0$$
 $$\phi _{Tc}(0,0) = x'_c(\frac{2\pi}{\sqrt{g'(0)}}, 0) = 0, \quad \psi  _{Tc}(0,0) = y'_c(\frac{2\pi}{\sqrt{g'(0)}}, 0) = \sqrt{g'(0)}$$ 
 
$$\phi _{cc}(0,0) = x_{cc}(\frac{2\pi}{\sqrt{g'(0)}}, 0), \quad \psi  _{cc}(0,0) = y_{cc}(\frac{2\pi}{\sqrt{g'(0)}}, 0) .$$

The derivatives satisfy the system
$$ 
 \cases{
 x'_{cc} = - {\sqrt{g'(0)}}y_{cc}  & \cr
 y'_{cc} =  \frac{1}{\sqrt{g'(0)}}[ [g'(x) - f'(x) y]x_{cc} + [g''(x) - f''(x) y](x_{c})^2 - f(x) y_{cc} - 2 f'(x) x_c y_c ]  & \cr}
$$

Setting now \ $c = 0$ \ in the preceding system, according to above intial conditions we get

$$ 
 \cases{
   x'_{cc} = - {\sqrt{g'(0)}}y_{cc}  & \cr
   y'_{cc} =  \frac{1}{\sqrt{g'(0)}}[ g'(0) x_{cc} + [g''(0) (\cos({\sqrt{g'(0)}} t))^2  - 2 f'(0) \cos({\sqrt{g'(0)}} t) \sin({\sqrt{g'(0)}} t) ]  & \cr}
$$
The solution of the latter system is

$$ 
 \cases{
   x'_{cc} =  \frac {1}{6 \sqrt{g'(0)}} [ - 3 g''(0) + g''(0) \cos({\sqrt{g'(0)}} t) + 4 f'(0) \sin({\sqrt{g'(0)}} t) + & \cr \qquad g''(0) \cos(2 {\sqrt{g'(0)}} t) - 2 f'(0) \sin(2 {\sqrt{g'(0)}} t) ] & \cr
   y'_{cc} =  \frac{1}{3 \sqrt{g'(0)}}[  [g''(0) \sin({\sqrt{g'(0)}} t)  - 2 f'(0) \cos({\sqrt{g'(0)}} t) + & \cr 
\qquad g''(0) \sin(2  {\sqrt{g'(0)}} t) + 2 f'(0) \cos({\sqrt{g'(0)}} t) ]  & \cr}
$$

We deduce from this the values
$$ \phi _{cc} (0,0) = \psi _{cc} (0,0) = 0$$ 
By similar method we establish the value of the third derivatives. We find in particular
$$\phi _{Tcc} (0,0) = x'_{cc} (\frac {2 \pi}{ \sqrt{g'(0)}}, 0) = 0, \qquad  \psi _{Tcc} (0,0) = y'_{cc} (\frac {2 \pi}{ \sqrt{g'(0)}}, 0) =     
\frac {g''(0)}{2 \sqrt{g'(0)}} $$

$$\phi _{TTc} (0,0) = x''_c (\frac {2 \pi}{ \sqrt{g'(0)}}, 0) = - g'(0), \qquad  \psi _{TTc} (0,0) = y''_c [\frac {2 \pi}{ \sqrt{g'(0)}}, 0) =  0    $$

$$\phi _{ccc} (0,0) = x_{ccc} (\frac {2 \pi}{ \sqrt{g'(0)}}, 0) = \frac {3 \pi}{2 \sqrt{g'(0)}} [ \frac {g''(0)}{2 g'(0)}\frac {f'(0)}{2\sqrt{g'(0)}} - \frac {f''(0)}{2\sqrt{g'(0)}} ]$$

$$\psi _{ccc} (0,0) = y_{ccc} (\frac {2 \pi}{ \sqrt{g'(0)}}, 0) = \frac { \pi}{2 \sqrt{g'(0)}} [ - \frac {f'^2(0)}{g'(0)}- \frac {10 g''^2(0)}{4 g'^2(0)} + 9 \frac {g^{(3)}(0)}{6 g'(0)} ]$$

From the above calculus, we may write the expansion of \ $\psi (T,c)$ \ near $(0,0)$
$$\phi (T,c) = - \frac{g'(0)}{2}(T - \frac {2 \pi}{ \sqrt{g'(0)}})^2 c +  \frac{1}{6}\phi _{ccc} (0,0) c^3 + o(c^3)$$

$$\psi (T,c) = (T - \frac {2 \pi}{ \sqrt{g'(0)}}) c + \frac {g''(0)}{2 \sqrt{g'(0)}}(T - \frac {2 \pi}{ \sqrt{g'(0)}}) c^2 + \frac{1}{6}\psi _{ccc} (0,0) c^3 + o(c^3)$$

Finally, return now to our problem. Solving \ $\psi (T(c),c) = 0$ \ for $T$ as function of $c$ 
$$ T(c) = \frac {2 \pi}{ \sqrt{g'(0)}} - \frac {1}{6} \psi _{ccc} (0,0) c^3 + o(c^3)$$
Recall the orbit is closed if \ $\Phi (c) = \phi (T(c), c) = 0.$\
In substituting \ $T = T(c)$\ we get \ $\Phi (c) =  \frac {1}{6} \phi _{ccc} (0,0) c^3 + o(c^3)$ 
Thus, we have the expansion of the period function for $c$ small
$$T(c) = \frac {2 \pi}{ \sqrt{g'(0)}} - \frac { \pi}{12 \sqrt{g'(0)}} [ - \frac {f'^2(0)}{g'(0)}- \frac {10 g''^2(0)}{4 g'^2(0)} + 9 \frac {g^{(3)}(0)}{6 g'(0)} ] c^2 + o(c^3)$$
which leads to the necessary condition and gives the theorem.

\bigskip

We obtain in particular this consequence

\bigskip

{\bf Corollary 2} \qquad {\it Let the conservative equation \ $x'' + g(x) = 0 $
\ and the Lienard equations\ $x'' \pm  f(x) x' + g(x) = 0$.\ Supppose the two functions  $f  \in C^1([a,b])$ , and $g \in C^3([a,b])$ , verifying\ $f(0) = g(0)  = 0, g'(0) > 0$\ and  let the origin $0$ be a center for these equations. \\ If the period function of \ $x'' + g(x) = 0 $
\ is increasing in a neighborhood of $0$,  then the period function of \ $x'' \pm  f(x) x' + g(x) = 0$\ is also increasing in a neighborhood of $0$.}

\bigskip
As a consequence of Theorem 1, we may deduce the following which also precises results of [S] .
\bigskip

{\bf Proposition 3} \qquad {\it  Let a function $f  \in C^2([a,b])$ , and $g \in C^3([a,b])$ , verifying\ $f(0) = g(0)  = 0, g'(0) > 0$\ and  let the origin $0$ be a center of (3.1) and in addition suppose that $g''(0) = f''(0) = 0$. Let \ $C(x) = \frac{1  }{g'(0) }g(x)  - \frac{1  }{g'(0) x^3} [\int_0^x s f(s)
ds]^2,$ \ we have in particular,
\begin{eqnarray}
{\it if}\  C(x) \ {\it is} \cases{ {\it \ strictly\ convex \ for}\  x\in J \ {\it  and}\  x < 0, &\cr
{\it \ strictly \  concave \ for}\ x\in J \ {\it and}\ x > 0, &\cr  }
\end{eqnarray}
then $T$ is increasing in a neighborhood of $0$;
 \begin{eqnarray}
{\it if}  \ C(x) \ {\it is} \cases{ {\it \ strictly \ convex \ for} \ x\in J \ {\rm and} \ x > 0, &\cr
{\it \ strictly \ concave \ for} x\in J \ {\rm and} \ x < 0, &\cr  }
\end{eqnarray}
then $T$ is decreasing in a neighborhood of $0$;
\begin{eqnarray}
  {\it  if} \  \frac{d^2  }{dx^2} C(x) \equiv 0 \ {\it for} \ x\in J,
\end{eqnarray}
then $T$ is constant in a neighborhood of $0$.} \\

{\bf Proof} \qquad After change in polar coordinates $(r,\theta
)$,  we obtain the following equivalent system to (3.4)
\begin{eqnarray}
\cases{
r' = \sqrt{g'(0)} r \cos \theta \sin \theta  - r \beta (r\cos \theta ) - \sin \theta  \ C(r \cos \theta ) & \cr
 \theta ' = - \frac{1}{\sqrt{g'(0)}} \cos ^2 \theta - \sin^2 \theta  - \frac{\cos \theta C(r \cos \theta ) - C'(0) (r \cos \theta)} {r}  & \cr }
\end{eqnarray}

We also observe that from system (3.1) we get $$r^2 \theta ' = r^2
\omega (r,\theta ) = - x C(x) - y^2,$$ here $T(r) = \int _{[0,2\pi
]} \frac{d\theta   }{\omega }$.\ By Theorem 1 of [S], it is
sufficient to prove for example that hypothesis (3.9) implies that
$\frac{\delta \omega (r,\theta )  }{\delta r} \leq 0$ for almost
all values $\theta  \in [0,2\pi ]$. 

A calculus gives $$ \frac{\delta \omega(r,\theta )  }{\delta r} = - \frac{\delta
}{\delta r} \frac{[\cos\theta (C(r\cos \theta )-C'(0) (r\cos\theta )]}{r} $$ 
$$= r \cos^2 \theta \ C'(r\cos \theta ) - \cos\theta \ C(r\cos \theta ).$$ 

We then obtain $$ \frac{\delta
\omega(r,\theta )  }{\delta r} = \frac{x C(x) -x^2 C'(x)  }{(x^2 +
y^2)^{3/2}} = x\frac{ C(x) -x C'(x)  }{(x^2 + y^2)^{3/2}}.$$
Remark obviously that since $(x C'(x) - C(x))' = x C''(x)$,  then
according to the hypothesis, the function $C \in C^3([a,b])$ and
condition $x C''(x) \leq 0$ for $x \in J$ which is equivalent to
hypothesis (3.9) implies $ \frac{\delta \omega(r,\theta )
}{\delta r} \leq 0$. \\ By the same way, we prove condition $x
C''(x) \geq 0$ for $x \in J$ which is equivalent to hypothesis
(3.6) implies $\frac{\delta \omega (r,\theta )  }{\delta r} > 0$.
\\ We then  prove the functions $\frac{\delta \omega   }{\delta r}
$ and $x C''(x)$ have a same sign. \\
In fact, we can that with a other manner.
A calculation gives $$C''(0) = \frac{g''(0)}{g'(0)} - \frac{2  }{3 g'(0)}f(0) f'(0) \quad {\it and} \quad C'''(0) = \frac{g'''(0)}{g'(0)} - \frac{2  }{3 g'(0)} f'(0)^2.$$
Then , we  necessarely have \ $C''(0) = 0$,\ since \ $xC''(x) \neq 0$ \ (if $x \neq 0$) implies  the monotonicity of the period.
 Thus, \ $x C''(x)$\ and  \ $g'''(0) - \frac{2  }{3} f'(0)^2$\ have the same sign in a neighborhood of $0$. This determines the monotonicity of the period function. \\ Also, condition $f''(0) = 0$ appears to be necessary by the following Lemma.
\\

\bigskip

{\bf Lemma 2} \qquad {\it Let  $f  \in C^2([a,b])$ , and $g \in C^3([a,b])$ , verifying\ $f(0) = g(0)  = 0, g'(0) > 0$\ and  let the Lienard equation  (3.1). Then, a necessary condition for the origin $0$ to be a center of (3.1) is }
$$ f'(0) g''(0) - 2 g'(0) f''(0) =  0 $$  

\bigskip
{\bf Proof}\qquad In the proof of Theorem 1, we have given the behavior of the period near the origin. \ $\Phi(c) \neq 0$ for small $c$ means there will be no any closed orbit in a neighborhood of  $0$. Indeed, $c + \Phi (c) $ corresponds to the position of return to the x-axis. Since,  we have seen 
$$\Phi (c) = \phi (T(c),c) =  - \frac{g'(0)}{2}(T - \frac {2 \pi}{ \sqrt{g'(0)}})^2 c +  \frac{1}{6}\phi _{ccc} (0,0) c^3 + o(c^3)$$ 
and according the expansion of $T(c)$, we get 
$$\Phi (c) = \frac{1}{6}\phi _{ccc} (0,0) c^3 + o(c^3)$$

Thus, a necessary condition for the origin to be a center is :

$$\phi _{ccc} (0,0)  = \frac {3 \pi}{2 \sqrt{g'(0)}} [ \frac {g''(0)}{2 g'(0)}\frac{f'(0)}{2\sqrt{g'(0)}} - \frac {f''(0)}{2\sqrt{g'(0)}} ] = 0$$

The lemma is then proved. 

\bigskip

Notice that by definition of \ $C(x),$ \ $f(x) \equiv 0$ \ implies \ $C(x) \equiv \frac{g(x)}{g'(0)}$\ and Corollary 2 reduces to the above Proposition 1. Remark also condition \ $g''(0) = 0$\ appears to be necessary to study the various monotonicity conditions of the period function. \\

The following which may be deduce from Theorem 1 is  particularly interesting
\bigskip

{\bf Corollary 3} \qquad {\it Let the conservative equation \ $x'' + g(x) = 0 $
\ and the Lienard equations\ $x'' \pm  f(x) x' + g(x) = 0$.\ Supppose the two functions  $f  \in C^1([a,b])$ , and $g \in C^3([a,b])$ , verifying\ $f(0) = g(0)  = 0, g'(0) > 0$\ and  let the origin $0$ be a center of these equations. \\ If \ $g''(0) \neq 0$ \ and the center $0$ of \ $x'' \pm  f(x) x' + g(x) = 0$\ is isochronous then the period function \ $T$ \ of \ $x'' + g(x) = 0 $
\ is strictly decreasing in a neighborhood of $0$.}

\bigskip
We prove that in considering Corrolary 2 , and on account of the fact that   $g$ is not a odd function. \\ According to this result, to establish the  existence of isochronous centers  for the Lienard sytem, it is necessary first to make sure that the period of the conservative associated system is decreasing at least in a neighborhood of $0$ .  One will be able to use for that the different criteria of monotonicity of the period function.\\

\bigskip
From Theorem 1 we may deduce other consequences.  \\
\bigskip

{\bf Corollary 4} \qquad {\it  Let a function $f  \in C^1([a,b])$ , and $g$ a smooth function, verifying\ $f(0) = g(0) , \ g'(0) > 0$\ and \ $g''(0) \neq 0, $ \  and  let the origin $0$ be a center of (3.1) . Then 
   necessary  conditions on the function \ $f$\ so that equation (3.1) has $T$  constant in a neighborhood of $0$ are
$$f'(0) = \pm \sqrt{3 g^{(3)}(0) - \frac{5}{g'(0)} g''^2(0)}, \qquad f''(0) = \pm \frac{g''(0)}{2 g'(0)} \sqrt{3 g^{(3)}(0) - \frac{5}{g'(0)} g''^2(0)} 
$$ }

\bigskip

Indeed, it suffices to remark by Corollary 3 the period function of the  equation\ $x'' + g(x) = 0$ \ has to be decreasing. Then, according to criteria of Schaaf it is necessary that \ $3 g^{(3)}(x) - \frac{5}{g'(x)} g''^2(x) > 0$\ in a neighborhood of $0$.\
This result is interesting since we do not need to suppose \ $g(x)$\  odd. 
In particular, if $g''(0) = 0$ we have  \ $\frac{d^2  }{dx^2} C(x) \equiv 0$\ implies   \ $C'''(0) = \frac{g^{(3)}(0)}{g'(0)} - \frac{2}{3 g'(0)} f'^2(0) = 0 $.\\

\bigskip

On the other hand, in considering the following function depending on $f$ and $g$
$$ \sigma (x) = 2x^2 \frac{1  }{g'(0) }f(x) \int _0^x sf(s) ds - 4\frac{1  }{g'(0) } [\int _0^x s f(s) ds ]^2 + x^3 g_n(x) - x^4 g'_n(x) $$
where $g_n(x) = \frac{1  }{g'(0) }g(x) -  x$.\   Then, for this $ \sigma (x) $  we have the following result which has been proved by [S] in the case where $g$ and $f$ are $C^1$   (Theorem A, Theorem 2 and Corollary 1)\\

\newpage

{\bf Proposition 2} \qquad {\it Let $f,g \in C^3(a,b), f(0) = g(0) = 0$ and $ g''(0) = f''(0) = 0$. The origin $0$ being a center of (3.1) . If \ $ x C(x) > 0$ \ in a punctured neighborhood $J$ of $0$, \\ then we have

(1) \quad if $\sigma (x) \leq 0$ for $x\in J$, then $T$ is decreasing in a neighborhood of $0$;

(2) \quad if $\sigma (x) \geq 0$ for $x\in J$, then $T$ is increasing in a neighborhood of $0$;

(3) \quad if $\sigma (x) \equiv 0$ for $x\in J$, then $T$ is constant in a neighborhood of $0$.\\ }

\bigskip

 Notice that while considering the assumption $xC(x) > 0$ in the case where $g$ is  $C^3$ and $f$ is  $C^1$, we have $g'(0) > 0$ and $C(x)$ is $C^3$ . 
Furthermore, $\sigma (x)$ and $C(x)$  are related $$\sigma (x) = - x^5 \frac{d}{dx} (\frac{C(x)}{x}.$$ This function may be written in a neighborhood of $0$,
 $$\sigma (x) = - 2 x^6 [g^{(3)}(x)  - \frac{2}{3} f'^2(x)] + ..... $$ That imposes an additional condition  \ $g''(0) = 0$ which turns out to be necessary since we have seen $C''(0) = g''(0)$. Moreover, it implies necessarely $f''(0) = 0$ by Lemma 2.\\

\subsection{Other class of equations} 

In [C2], Chicone consider the class of differential equations of the form $$x'' + F(x') + G(x) = 0$$ with Dirichlet or Neumann boundary values. These functions  verify   $F(0) = G(0) = 0$ . This is equivalent to the system 

\begin{eqnarray}
   \cases{
   x' = - y  & \cr
   y' =   x - x \tilde{g}(x) - y \tilde{f}(y)  & \cr }
\end{eqnarray}
   
In the stantard Neuman situation the functions $\tilde{f}$ and $\tilde{g}$ are such that $\tilde{f},\tilde{g} \in C^2([-a,a])$ and satisfy the conditions \ {\bf (C)}
$$\quad \tilde{f}(-s) = - \tilde{f}(s), \qquad \tilde{g}(-s) = - \tilde{g}(s), \quad {\it for } \quad s \in [-a,a]$$ and $$\tilde{f}'(s) \geq 0, \quad \tilde{g}'(s) \geq 0, \quad \tilde{f}''(s) \geq 0, \quad \tilde{g}''(s) \geq 0, \quad {\it for} \quad s \in [0,a].$$ 

The trajectories of this system are symmetric with respect to the x-axis. So, it has obviously a center at the origin of the phase plane. Chicone proved that under the above conditions, the period function is monotone increasing. \\
The following class of Raleigh equations with linear restoring term may have a monotone increasing period  in a  neighborhood of the center $0$. 
$$ (R) \qquad x'' + F(x') + x  = 0 .$$

Without supposing  hypotheses above, we may prove an analogous result which also improves Corollary 9 of [S]. We only need for that to suppose $F(x)$  is an even function with \ $ F''(0)  \neq 0, $\ instead of conditions $\bf (C)$. \\

 \bigskip

{\bf Corollary 5} \qquad {\it Let $h$ be an analytic even function such that \ $F(0) = 0,$ \ and \ $F''(0) \neq 0$.\ The origin $0$ be a center of $(R)$. Then the period function \ $T$ \ of equation \ $R$\ is monotone increasing in a  neighborhood of $0$. }
\bigskip

{\bf Proof} \qquad  The system 

\begin{eqnarray}
   \cases{
   x' = - y  & \cr
   y' =   x  +  F(y)  & \cr }
\end{eqnarray}

which is equivalent to (R), has a unique singular point at the origin. By exchanging the variables and multiplying by $-1$ we obtain the equivalent system 

\begin{eqnarray}
   \cases{
   x' = - y  -  F(x)  & \cr
   y' =  x    & \cr }
\end{eqnarray}

which is a lienard system with \ $g(x) = x$,\ and \ $f(x) = \pm F'(x)$.\  Moreover, the origin is a center implies \ $f(x)$ \ is odd (Corollary 7 of [S]). The conclusion holds by Theorem 1 since \ $ - \frac{2}{3} f'^2(0)g'(0) < 0$.

\bigskip

\subsection{ General remarks}

The results that we propose (Theorem 1)  presents the advantage of to be natural and easy to apply.  The monotonicity condition as well as the isochronicity case can verify itself more easily.  So much as the function $C(x)$  reduces itself to the function $g(x)$  when $f(x) \equiv 0$ in a neighborhood of $0$.  This puts in evidence the parallel between our result and the one of Opial. More precisely, a calculation gives $$C''(0) = g''(0) - \frac{2  }{3}f(0) f'(0) \quad {\it and} \quad C'''(0) = g'''(0) - \frac{2  }{3} f'(0)^2.$$ Then we have necessarely $C''(0) = 0$. \\ 
We may remark that  condition $C'''(x) < 0$ for $x \in (a,b)$ implies (2.7)  and $C'''(x) > 0$ implies (2.8) (this requires  necessarily  $C'''(0) < 0)$ when $g$ is $C^3$ and $g'''(0) < 0$). \\ It  naturally seems  that the function $C(x)$ plays a same role for the system (3.3) as the function $g(x)$ for the conservative system. Indeed, if we take $f(x) \equiv 0$, conditions (3.7) and (3.8) reduce to the conditions (2.2) and (2.3) respectively of the Proposition 1. Notice that Rothe condition for the monotonicity of the period function
$$H_{4}(x)=x[3g^{\prime}(x)^{2}-g(x)g^{\prime\prime}(x)-(3\frac{g^{\prime
}(0)^{2}}{g^{\prime\prime}(0)})g^{\prime\prime}(x)]\geq0 (\leq 0) \ ,
\  {\it for } \qquad {x\in(a,b),}
$$
is more general than (2.2) and (2.3) (see Remark 1 and [R]). \\ We may expect that the results of the above Theorem 1 can be generalize. In the sense where \ $C(x)$\ is replaced by another more general appropriated  function let \ $D(x)$,\ which itself can be reduce to the Rothe function  $$(3\frac{g^{\prime
}(0)^{2}}{g^{\prime\prime}(0)})g^{\prime\prime}(x) - 3g^{\prime}(x)^{2} + g(x)g^{\prime\prime}(x)$$
for the conservative case $f(x) \equiv 0.$\ This function have to be  such that $$ D^{(3)}(0) =  g^{(3)}(0) - \frac{5}{3g'(0)} g''^2(0) - \frac{2}{3} f'^2(0)$$
One separates thus the necessary condition $g''(0) = 0.$ \\ We may also expect that the sign of the following function $$g'(x) g^{(3)}(x) - \frac{5}{3} g''^2(x) - \frac{2}{3} f'^2(x) g'(x) $$ determines the global monotonicity of the period function $T$ of the Lienard system.  In the conservative case, the last function reduces to   $g'(x) g^{(3)}(x) - \frac{5}{3} g''^2(x)$  which intervenes in the Schaaf condition ( see condition $(\mathcal{C}_3)$ in Section 2 above). Recall that the last one is less inclusive than the Rothe condition.\\ Moreover, according Corollary 4 to determine the isochronous centers  at the origin for the Lienard system (other than those determined by [S] ) we have to insure that the associated conservative sytem has a decreasing period function. \\

\newpage

\begin{center}
\Large
REFERENCES\\
\end{center}

\bigskip

[C] \qquad C.Chicone \quad {\it The monotonicity of the period function for planar Hamiltonian vector fields} \quad J. Diff Eq., vol 69, p. 310-321, (1987).\\

[C2] \qquad C.Chicone \quad {\it Geometric methods for two-point nonlinear boundary value problems} \quad J. Diff Eq., vol 72, p. 360-407, (1988).\\

[Ch-C] \qquad R.Chouikha and F. Cuvelier \quad {\it Remarks on some monotonicity conditions for the period function}\quad  Applic. Math., vol 26, 3, p. 243-252, (1999).\\

[C-W] \qquad S.N. Chow and D. Wang \quad {\it On the monotonicity of the period function of some second order equations}\quad Casopis Pest. Mat. 111, p. 14-25, (1986).\\

[C-D] \qquad CJ Christopher and J. Devlin \quad {\it Isochronous centers in planar polynomial systems} \quad SIAM J. Math. Anal., vol 28, p.162-177, (1997).\\

[L] \qquad W.S. Loud \quad {\it The behavior of the period of solutions ofcertain plane autonomous systems near centers  }\quad Contr. Differential Equations, 3, p. 21-36, (1964). \\

[O] \qquad Z. Opial \quad {\it Sur les periodes des solutions de l'equation differentielle \ $x'' + g(x) = 0$} \quad Ann. Polon. Math., 10, p. 49-72, (1961).\\

[R] \qquad F. Rothe \quad {\it Remarks on periods of planar Hamiltonian systems.} \quad SIAM J. Math. Anal., 24, p.129-154, (1993).\\

[S] \qquad M. Sabatini \quad {\it On the period function of Lienard systems} \quad J. of Diff. Eq., 152, p. 467-487, (1999).

[Sc] \qquad R.Schaaf \quad {\it A class of Hamiltonian systems with increasing periods} \quad J. Reine Angew. Math., 363, p. 96-109, (1985).\\

\end{document}